\documentclass[11pt]{article}
\usepackage {jr,myfigures,enumerate}%
\usepackage {apajr}
\def\eq{}

\title {Hierarchical selection of variables in sparse
high-dimensional regression}

\author{P. J. Bickel\\Department of Statistics\\University of California at Berkeley
\and Y. Ritov\\ Department of Statistics\\The Hebrew University of
Jerusalem \and A.B. Tsybakov\\Laboratoire de
Statistique, CREST, Timbre J340 \\
3, av.Pierre Larousse, 92240 Malakoff cedex, France \\
and Laboratoire de Probablit\'es et Mod\`eles
Al\'eatoires\\Univerist\'e Pierre et Marie Curie}

\def\ms{\scm\scs}

\begin{document}
\maketitle
\begin{abstract}
 We study a regression model with a huge number of interacting
 variables. We consider a specific approximation of the regression
 function under two assumptions: (i) there exists a sparse representation
 of the regression function in a suggested basis, (ii) there are
no interactions outside of the set of the corresponding main
effects. We suggest an hierarchical randomized search procedure for
selection of variables and of their interactions.   We show that given an initial estimator, an estimator with a similar prediction loss but with a smaller number of non-zero coordinates can be found.

\end{abstract}

\section{Introduction}
Suppose that we observe $(Y_i,\mbX_i)$, $i=1,\dots,n$, an \iid
sample from the joint distribution of $(Y,\mbX)$, where $Y\in\scr$,
and $\mbX=(X_1,\dots,X_d)\in\scx_1\times\cdots\times\scx_d=\scx$,
with $\scx_j$ being some subsets of finite-dimensional Euclidean
spaces. Our purpose is to estimate the regression function
$f(\mbX)=E(Y|\mbX)$ nonparametrically by constructing a suitable
parametric approximation of this function, with data-dependent
values of the parameters. We consider the situation where $n$ is
large, or even very large and the dimension $d$ is also large.
Without any assumptions, the problem is cursed by its dimensionality
even when $\scx_j=\scr$ for all $j$. For example, a histogram
approximation has $p=3^{20}>10^9$ parameters when the number of
variables is $d=20$, and the range of each is divided into the meager
number of three histogram bins.

It is common now to consider models where the number of parameters
$p$ is much larger than the sample size $n$. The idea is that the
effective dimension is defined not by the number of potential
parameters $p$ but by the (unknown) number of non-zero parameters
that can be much smaller than $n$. Methods like thresholding in
white noise model, cf. \cite{abdj06} or \cite{g02}, LASSO, LARS or
Dantzig selector in regression, cf, \cite{t96}, \cite{cds01},
\cite{ehjt04}, \cite{ct07}, are used, and it is proved that if the
vector of estimated parameters is sparse (i.e., the number of
non-zero parameters is relatively small) then the model can be
estimated with reasonable accuracy, cf.
\cite{btw07,btw06b,ct07,kf00,gr04,mb06,my06,zh06,zy05}. A direct
selection of a small number of non-zero variables is relatively
simple for the white noise model. There, each variable is processed
separately, and the parameters can be ordered according to the
likelihood that they are non-zero. The situation is more complicated
in regression problems. Methods like LASSO and LARS yield
numerically efficient ways to construct a sparse model, cf.
\cite{jn00,n00,os00a,os00b,ehjt04,t05}. However, they have their
limits, and are not \emph{numerically} feasible with too many
parameters, as for instance in the simple example considered above.

Our aim is to propose a procedure that can work efficiently
in such situations. We now outline its general scheme. 
Consider a collection of functions
$(\psi_{i,j})_{i=1,\dots,d,\,j=0,1,\dots,L}$ where
 $\psi_{i,j}:\scx_i\to\scr$. For example, for fixed
$i$ this can be a part of a basis $(\psi_{i,j})_{j=0,1,\dots}$ for
$L_2(\scx_i)$. For simplicity, we take the same number $L$ of basis
functions for each variable. We assume that $\psi_{i,0}\equiv 1$.
Consider an approximation $f_\beta$ of regression function $f$ given
by: \eqsplit{
    f_\beta(\mbX)=
    \sum_{\mbj\in \{0,1,\dots,L\}^d}
    \beta_{\mbj}\prod_{i=1}^d\psi_{i,j_i}
    (X_{i})
  }
 where $\mbj = (j_1,\dots,j_d)$ and $\beta_{\mbj}$ are unknown
 coefficients.   Note that $f_\beta$ is nothing but a specific model with
interactions between variables, such that all the interactions are
expressed by products of functions of a single variable. In fact,
since $\psi_{i,0}\equiv 1$, the multi-indices $\mbj$ with only one
non-zero coefficient yield all the functions of a single variable,
those with only two non-zero coefficients yield all the products of
two such functions, etc. Clearly, this covers the above histogram
example, wavelet approximations and others.

The number of coefficients $\beta_{\mbj}$ in the model is $(L+1)^d$.
The LASSO type estimator can deal with a large number of potential
coefficients which grows exponentially in $n$. So, theoretically, we
could throw all the factors into the LASSO algorithm and find a
solution. But $p \sim L^d$ is typically a huge number. Although in
the theory LASSO can handle that many variables, in practice, it
becomes numerically infeasible. Therefore, a systematic search is
needed.

Since there is no way to know in advance which factors are
significant, we suggest a hierarchical selection: we build the model
in a tree fashion.  At each step of the iteration we apply a LASSO
type algorithm to a collection of candidate functions, where we
start with all functions of a single variable. Then, from the model
selected by this algorithm we extract a sub-model which includes
only $K$ functions, for some predefined $K$. The next step of the
iteration starts with the same candidate functions as its
predecessor plus all the interactions between the $K$ functions
selected at the previous step.

Formally we consider the following hierarchical model selection
method. For a set of functions $\scf$ with cardinality $|\scf|\ge
K$, let $\ms_K$ be some procedure to select $K$ functions out of
$\scf$. We denote by $\ms_K(\scf)$ the selected subset of $\scf$,
$|\ms_K(\scf)|=K$. Also, for a function $f:\scx\to\scr$, let
$\bbn(f)$ be the minimal set of indices  such that $f$ is a function
of $(X_i)_{i\in\bbn(f)}$ only. The procedure is defined as follows.
\begin{enumerate}
\item[(i)]
Set $\scf_0=\union_{i=1}^d\{\psi_{i,1},\dots,\psi_{i,L}\}$.
\item[(ii)] For $m=1,2,\dots$,
let $$\scf_m = \scf_{m-1}\union\{fg:\; f,g\in\scm\scs_K(\scf_{m-1}),
\bbn(f)\inter\bbn(g)=\emptyset \}.$$
\item[(iii)] Continue until convergence is declared.
The output of the algorithm is the set of functions
$\scm\scs_K(\scf_m)$ for some $m$.
\end{enumerate}

This search procedure is valid under the dictum of no interaction
outside of the set of the corresponding main effects: a term is
included only if it is a function of one variable or it is a product
of two other included terms. If this is not a valid assumption one
can enrich the search at each step to cover  all the
coefficients $\beta_{\mbj}$ of the model. However, this would be
cumbersome.

Note  that $|\scf_m|\leq K^2+|\scf_{m-1}| \leq m K^2 + |\scf_0| = m
K^2 + Ld$. Thus, the set $\scf_m$ is not excessively large. At every
step of the procedure we keep for selection all the functions of a
single variable, along with not too many interaction terms. In other
words, functions of a single variable are treated as privileged
contributors. On the contrary, interactions are considered with a
suspicion increasing as their multiplicity grows: they cannot be
candidates for inclusion unless their ``ancestors" were included at
all the previous steps.

The final number of selected effects is $K$ by construction. We
should choose $K$ to be much smaller than $n$ if we want to fit our
final model in the framework of the classical regression theory.

One can split the sample in two parts and do model selection and
estimation separately. Theoretically, the rate of convergence of the
LASSO type procedures suffers very little when the procedures are
applied only to a sub-sample of the observations, as long as the
sub-sample size $n_{MS}$ used for model selection is such that
$n_{MS}/n$ converges slowly to 0.  We can therefore, first use a
sub-sample of size $n_{MS}$ to select, according to (i)--(iii), a
set of $K$ terms that we include in the model. The second stage will
use the rest of the sample and estimate via, e.g., standard
least-square method the regression coefficients of the $K$ selected
terms.

This paper has two goals. The first one, as described already, is
suggesting a method to build highly complex models in a hierarchial
fashion.  The second purpose is arguing that a reasonable way to do
model selection is a two stage procedure. The first stage can be
based on the LASSO, which is an efficient way to obtain sparse
representation of a regression model. We argue, however, by a way of
example in Section \ref{sec:exam}, that using solely the LASSO can
be an non-optimal procedure for model selection. Therefore, in
Section \ref{rs} we introduce the second stage of selection, such
that a model of a desired size is obtained at the end. At this stage
we suggest to use either randomized methods or the standard backward
procedure. We prove prediction error bounds for two randomized
methods of pruning the result of the LASSO stage. Finally, in
Section \ref{sec:examples} we consider two examples that combine the
ideas presented in this paper.

\section{Model selection: an example}\label{sec:exam}

The above hierarchical method depends on a model selection procedure
$\ms_K$ that we need to determine. For high-dimensional case that we
are dealing with, LASSO is known to be an efficient model selection
tool: it is shown that under general conditions the set of non-zero
coefficients of LASSO estimator coincides with the true set of
non-zero coefficients in linear regression, with probability
converging to 1 as $n\to\infty$ (see, e.g., \cite{mb06,zy05}).
However, these results depend on strong assumptions that essentially
role off anything close to multicolinearity. These conditions are
often violated in practice when there are many variables
representing a plentitude of highly related one to another
demographic and physical measurements of the same subject. They are
also violated in a common statistical learning setup where the
variables of the analysis are values of different functions of one
real variable (e.g., different step functions).
 Note
that for our procedure we do not need to retain all the non-zero
coefficients but just to extract the $K$ ``most important" ones. To
achieve this, we first tried to tune the LASSO in some natural way.
However, this approach failed.

We start with an example. We use this example to argue that although
the LASSO does select a small model (i.e., typically many of the
coordinates of the LASSO estimator are 0), it does a poor job in
selecting the relevant variables. A naive approach for model
selection when the constraint applies to the number of non-zero
coefficients, is to relax the LASSO algorithm until it yields a
solution with the right number of variables. We believe that this is
a wrong approach. The LASSO is geared for $L_1$ constraints and not
for $L_0$ ones. We suggest another procedure in which we run the
LASSO until it yields a model more complex than wished, but not too
complex, so that a standard model selection technique like backward
selection can be used. This was the method considered in \cite{gr04}
to argue that there are model selection methods which are persistent
under general conditions.

We first recall the basic definition of LASSO. Consider the linear
regression model $$\mby=\mbZ\beta_0+\varepsilon$$ where
$\mby=(Y_1,\dots,Y_n)' \in \scr^n$ is the vector of observed
responses, $\mbZ\in\scr^{n\times p}$ is the design matrix,
$\beta_0\in\scr^p $ is an unknown parameter and $\varepsilon=
(\xi_1,\dots,\xi_n)'\in \scr^n$ is a noise. The LASSO estimator
$\hat\beta_{L}$ of $\beta_0$ is defined as a solution of the
minimization problem
\begin{equation}\label{las1}
\min_{\beta: \,\|\beta\|_1\le T}\|\mby-\mbZ\beta\|^2
\end{equation}
where $T>0$ is a tuning parameter, $\|\beta\|_1$ is the
$\ell_1$-norm of $\beta$ and $\|\cdot\|$ is the empirical norm
associated to the sample of size $n$: $$\|\mby\|^2 =
n^{-1}\sum_{i=1}^n Y_i^2.$$ This is the formulation of the LASSO as
given in \cite{t96}. Another formulation, given below in
\eqref{las2}, is that of minimization of the sum of squares with
$L_1$ penalty. Clearly, \eqref{las1} is equivalent to \eqref{las2}
with some constant $r$ dependent on $T$ and on the data, by the
Lagrange argument. The standard LARS-like algorithm of
\cite{ehjt04}, which is the algorithm we used, is based on gradual
relaxation of the constraint $T$ of equation \eqref{las1}, and
solves therefore simultaneously both problems. The focus of this
paper is the selection of a model of a given size. Hence we apply
the LARS algorithm until we get for the first time a model of a
prescribed size.

\begin{example}
We consider a linear regression model with $100$ \iid observations
of $(Y,Z_1,\dots, Z_{150})$ where the predictors
$(Z_1,\dots,Z_{150})$ are \iid standard normal, the response
variable is $Y=\summ j1{150} \beta_j Z_j+\xi= \summ
j1{10}\frac{10}{25+j2}Z_j+\xi$, and the measurement error is
$\xi\dist N(0,\sig^2)$, $\sig=0.1$.

Note that we have more variables than observations but most of the
$\beta_j$ are zero.

Figure \ref{Fig:SelVarFun2}\hspace{-0.35em}a  presents the
regularization path, i.e. the values of the coefficients of
$\hat\beta_{L}$ as a function of $T$ in \eqref{las1}.
 The vertical dashed lines indicate the values of the
$T$ for which the number of non-zero coefficients of $\hat\beta_{L}$
is for the time larger than the mark value (multiple values of 5).
The legend on the right gives the value of the 20 coefficients with
the highest values (sorted by the absolute value of the
coefficient).

 Figure \ref{Fig:SelVarFun2}\hspace{-0.35em}b presents a similar
situation. In fact, the only difference is that the correlation
between any two $Z_i$'s is now 0.5. Again, the 10 most important
variables are those with non-zero true values.

 \twofiguresV[2.775]{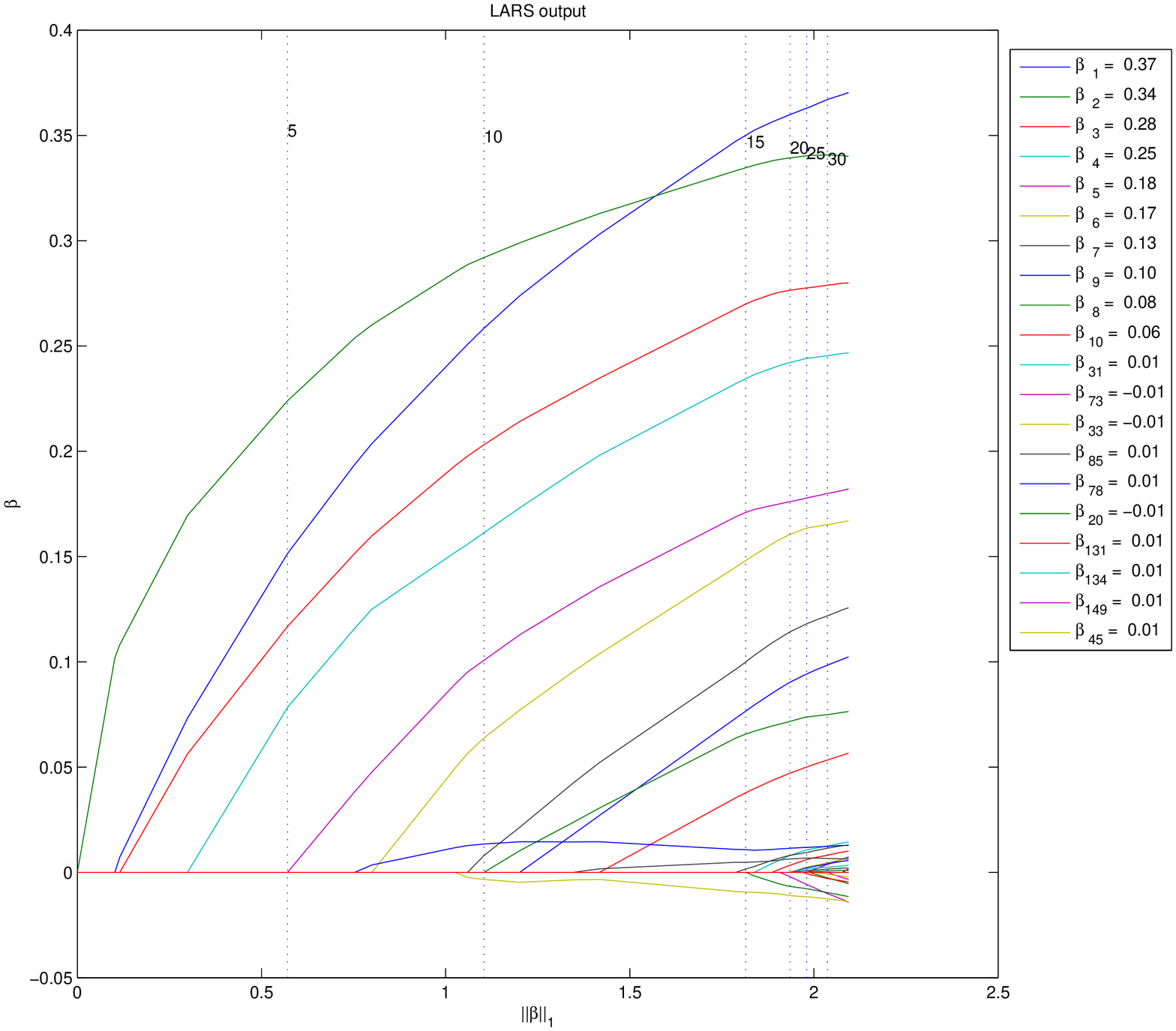}{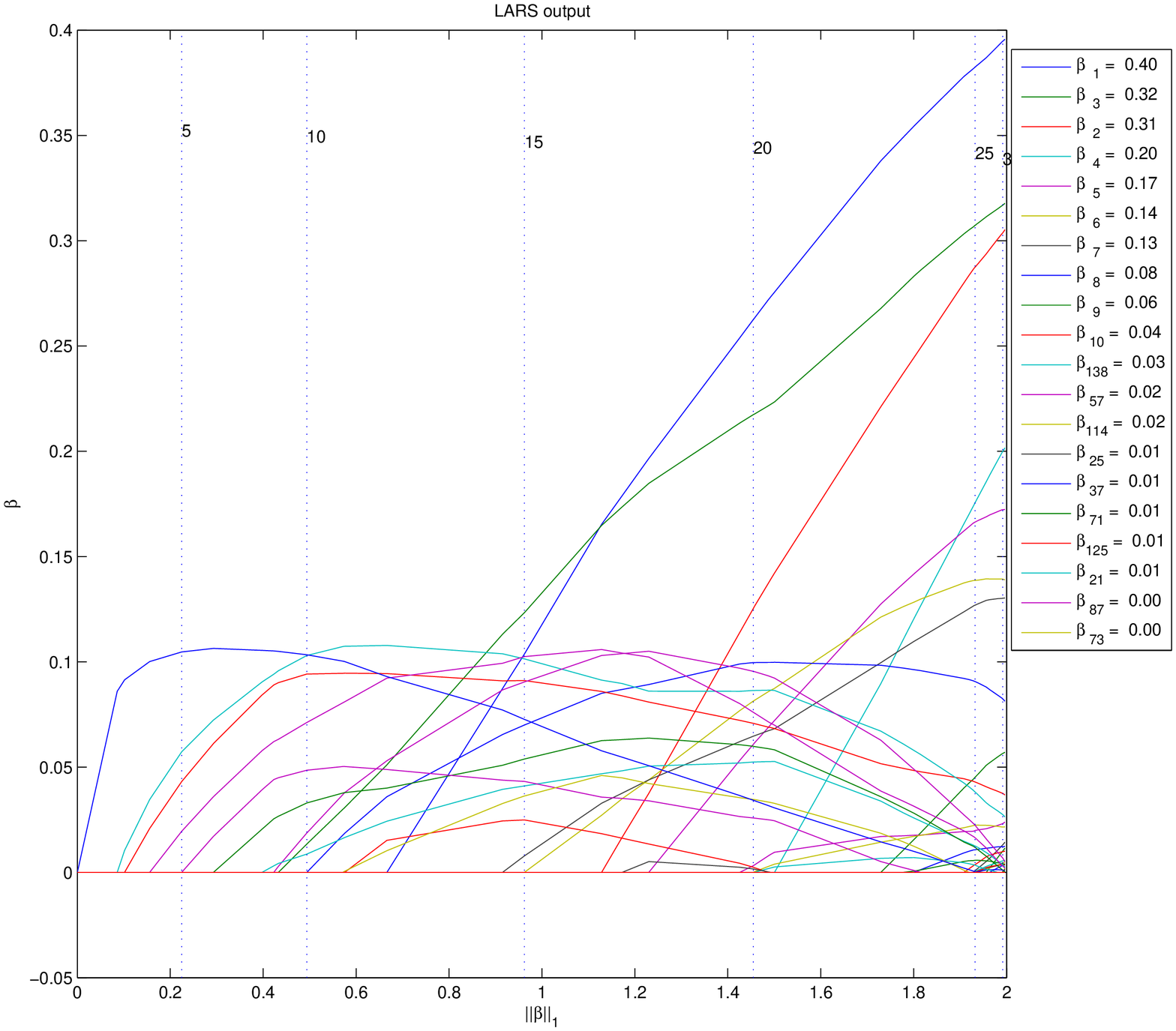}
    {Selecting variables. Coefficients vs. $L_1$ }{Fig:SelVarFun2}

\noindent

Suppose we knew in advance that there are exactly 10 non-zero
coefficients. It could be assumed that LASSO can be used, stopped
when it first finds 10 non-zero coefficients (this corresponds to
$T\approx 0.5$ in Figure \ref{Fig:SelVarFun2}\hspace{-0.35em}b).
However, if that was the algorithm, then only three coefficients
with non-zero true value, $\beta_3$, $\beta_8$, and $\beta_{10}$,
were included together with some 7 unrelated variables. For
$T\approx 2$ the 10 largest coefficients do correspond to the 10
relevant variables, but along with them many unrelated variables are
still selected (8 variables in Figure
\ref{Fig:SelVarFun2}\hspace{-0.35em}b), and moreover this particular
choice of $T$ cannot be known in advance if we deal with real data.
\end{example}

\section{Randomized selection}\label{rs}

The approach to design the model selector $\ms_K$ that we believe
should be used is the one applied in the examples of Section
\ref{sec:examples}. It acts as follows: run the LASSO for a large
model which is strictly larger than the model we want to consider,
yet small enough so that standard methods for selecting a good
subset of the variables can be implemented. Then run one of such
methods, with given subset size $K$: in the examples of Section
\ref{sec:examples} we use the standard backward selection procedure.
We do not have a mathematical proof which is directly relevant to
such a method. We can prove, however, the validity of an inferior
backward method which is based on random selection (with appropriate
weights) of the variable to be dropped at each stage. We bound the
increase in the sum of squares of the randomized method. The same
bounds are applied necessarily to the standard backward selection.

 Suppose that we have
an arbitrary estimator $\ti\beta$ with values in $\scr^p$, not
necessarily the LASSO estimator. We may think, for example, of any
estimator of parameter $\beta_0$ in the linear model of Section
\ref{sec:exam}, but our argument is not restricted to that case. We
now propose a randomized estimator $\widehat\beta$ such that:
\begin{itemize}
\item[(A)] the prediction risk of $\widehat\beta$ is on the average not too far from
that of $\ti\beta$,
\item[(B)] $\widehat\beta$ has at most $K$ non-zero components,
\item[(C)] large in absolute value components of
$\widehat\beta$ coincide with those of $\ti\beta$.
\end{itemize}

\medskip

{\it Definition of the randomization distribution.} Let $\sci$ be
the set of non-zero coordinates of the vector
$\ti\beta=(\ti\beta_1,\dots,\ti\beta_p)$. We suppose that its
cardinality $\ti K= |\sci|\ge 2$. Introduce the values
$$p_i =
\min\{1,c(\ti K-1)|\ti\beta_i|/\|\ti\beta\|_1\}, \quad i\in\sci,
$$
where $c\geq1$ is a solution of $\sum_{i\in\sci} p_i=\ti K-1$. Such
$c$ exists since the function $$t\mapsto \bar p_i(t) \equiv
\min\{1,t(\ti K-1)|\ti\beta_i|/\|\ti\beta\|_1\}$$
is continuous and
non-decreasing, $\lim_{t\to\infty} \sum_{i\in\sci} \bar p_i(t) = \ti
K$ and $\sum_{i\in\sci} \bar p_i(1) \le \ti K-1$. From
$\sum_{i\in\sci} p_i=\ti K-1$ we get
\begin{equation}\label{eq:sum1} \sum_{i\in\sci}(1-p_i)=1,
\end{equation}
so that the collection $\{1-p_i\}_{i\in\sci}$ defines a probability
distribution on $\sci$ that we denote by ${\rm P}^*$. Note that
there exists a $p_i$ not equal to 1 (otherwise we have
$\sum_{i\in\sci} p_i=\ti K$), in particular, we have always $p_i<1$
for the index $i$ that corresponds to the smallest in absolute value
$\ti\beta_i$. On the other hand, $p_i>0$ since $\ti\beta_i\neq 0$
for $i\in \sci$. Therefore, $0<p_i<1$ for at least two indices $i$
corresponding to the two smallest in absolute values coordinates of
$\ti\beta$.

\medskip

{\it Definition of the randomized selection procedure.} Choose $i^*$
from \sci at random according to distribution ${\rm P}^*$: ${\rm
P}^*(i^*=i)=1-p_i$, $i\in\sci$. We suppose that the random variable
$i^*$ is independent of the data ${\bf y}$. Define a randomized
estimator $\beta^*=(\beta^*_1,\dots,\beta^*_p)$ where
$\beta^*_{i^*}=0$, $\beta^*_i=\ti\beta_i/p_i$ for
$i\in\sci\setminus\{i^*\}$, and $\beta^*_i=0$ for $i\not\in \sci$.
In words, we set to zero one coordinate of $\ti\beta$ chosen at
random, and the other coordinates are either increased in absolute
value or left intact. We will see that on the average we do not
loose much in prediction quality by dropping a single coordinate in
this way.

We then perform the same randomization process taking $\beta^*$ as
initial estimator and taking randomization independently of the one
used on the first step. We thus drop one more coordinate, etc.
Continuing iteratively after $\ti K-K$ steps we are left with the
estimator which has exactly the prescribed number $K$ of non-zero
coordinates. We denote this final randomized estimator by
$\widehat{\beta}$. This is the one we are interested in.

\medskip

Denote by $\E^*$ the expectation operator with respect to the
overall randomization measure which is the product of randomization
measures over the $\ti K-K$ iterations.
\begin{theorem}\label{th1}
Let $\mbZ\in \scr^{n\times p}$ be a given matrix. Suppose that the
diagonal elements of the corresponding Gram matrix $\mbZ'\mbZ/n$ are
equal to 1, and let $\ti\beta$ be any estimator with $\ti K\ge 3$
non-zero components. Then the randomized estimator $\widehat\beta$
having at most $K<\ti K$ non-zero coordinates has the following
properties.
\begin{itemize}
\item[(i)] For any vector ${\bf f}\in\scr^n$,
$$\E^*\|{\bf f} - \mbZ\widehat\beta\|^2 \le
\|{\bf f} - \mbZ\ti\beta\|^2+\|\ti\beta\|_1^2
   \left(\frac{1}{K-1} - \frac{1}{\ti K-1}\right)\, .$$
\item[(ii)]
Let  $\ti\beta_{(j)}$ be the coordinates of $\ti\beta$ ordered by
absolute value: $|\ti\beta_{(1)}|\geq |\ti\beta_{(2)}|\geq\dots\ge
|\ti\beta_{(p)}|$. Suppose that
$|\ti\beta_{(k)}|>\|\ti\beta\|_1/(\ti K-1)$ for some $k$. Then the
estimator $\widehat\beta$ coincides with $\ti\beta$ in the $k$
largest coordinates: $\widehat\beta_{(j)}=\ti\beta_{(j)}$,
$j=1,\dots,k$.
\item[(iii)]
Suppose that $|\ti\beta_{(k+1)}|=0$ and
$|\ti\beta_{(k)}|>\|\ti\beta\|_1/(\ti K-1)$ for some $k$. Then
$\widehat\beta$ keeps all the non-zero coordinates of $\ti\beta$.
\end{itemize}\end{theorem}
\begin{proof}
 It is easy to see that
$\E^*(\beta^*_i)=\ti\beta_i$ for all $i$ and, for any vector ${\bf
f}\in\scr^n$, \eqsplit[ranvar]{
    \E^* \|{\bf f} - \mbZ\beta^*\|^2
    &=\|{\bf f}-\mbZ\ti\beta\|^2 + {1\over n}\trace(\mbZ'\mbZ\Sig^*)
    \\
    &= \|{\bf f}-\mbZ\ti\beta\|^2 + {1\over
    n}\sum_{i=1}^n\mbz_i'\Sig^*\mbz_i\\
    &\leq \|{\bf f}-\mbZ\ti\beta\|^2 + \sum_{j=1}^p\ti\beta_j^2
    \frac{1-p_j}{p_j}
  }
where $\mbz_i$ are the rows of matrix $\mbZ$ and
$\Sig^*=\E^*[(\beta^*-\ti\beta)(\beta^*-\ti\beta)']$ is the
randomization covariance matrix. We used here that $\Sig^*$
 is of the form
$$
\Sig^*={\rm diag}\left(\ti\beta_j^2
    \frac{1-p_j}{p_j}\right) - (B\ti
    \beta)(B\ti\beta)'\quad \text{with} \quad
    B={\rm diag}\left(\frac{1-p_i}{p_i}\right),
$$
and the diagonal elements of $\mbZ'\mbZ/n$ are equal to 1, by
assumption of the theorem.

Recall that $c\ge1$, and therefore $|\ti\beta_j|\geq
\|\ti\beta\|_1/(\ti K-1)$ implies $p_j=1$. Hence,
\eqsplit[vhat]{\eq
    \sum_{j\in\sci} \ti\beta_j^2
    \frac{1-p_j}{p_j} & = \sum_{0<|\ti\beta_j|<\|\ti\beta\|_1/(\ti K-1)}
    \ti\beta_j^2
    \frac{1-p_j}{p_j}
    \\
    &\le\frac{\|\ti\beta\|_1}{c(\ti K-1)}
    \sum_{0<|\ti\beta_j|<\|\ti\beta\|_1/(\ti K-1)}
    |\ti\beta_j|(1-p_j)
    \\
    &\leq \frac{\|\ti\beta\|_1^2}{(\ti K-1)^2}\sum_{j\in\sci}(1-p_j)
    \\
    &= \frac{\|\ti\beta\|_1^2}{(\ti K-1)^2}
    }
where we used \eqref{eq:sum1}. Thus, the randomized estimator
$\beta^*$ with at most $\ti K-1$ non-zero components satisfies
 \eqsplit[changedLS]{
    \E^* \|\mbf - \mbZ\beta^*\|^2
    &\leq \|\mbf-\mbZ\ti\beta\|^2 +
    \frac{\|\ti\beta\|_1^2}{(\ti K-1)^{2}}\,.
  }
Note also that $\beta^*$ has the same $\ell_1$ norm as the initial
estimator $\ti\beta$: \eqsplit[newL1]{
    \|\beta^*\|_1 &= \|\ti\beta\|_1
   }
In fact, the definition of $\beta^*$ yields
 \eqsplit{
    \|\beta^*\|_1 - \|\ti\beta\|_1 &= \left(\sum_{j\in\sci}
    \frac{|\ti\beta_j|}{p_j} -\frac{|\ti\beta_{i^*}|}{p_{i^*}}\right)
    - \sum_{j\in\sci}  |\ti\beta_j|
    \\
    &=  \frac{\|\ti\beta\|_1}{c(\ti K-1)} \sum_{p_j<1}
    \left(1 - c(\ti K-1)\frac{|\ti\beta_j|}{\|\ti\beta\|_1}\right)
    -  \frac{\|\ti\beta\|_1}{c(\ti K-1)}
    \\
    &= \frac{\|\ti\beta\|_1}{c(\ti K-1)}\sum_{j\in\sci}
    \Bigl(1 - p_j\Bigr)  - \frac{\|\ti\beta\|_1}{c(\ti K-1)}
    \\
    &= 0,
    }
in view of \ref{eq:sum1}.

Using \eqref{changedLS} and \eqref{newL1} and continuing by
induction we get that the final randomized estimator $\widehat\beta$
satisfies \eqsplit{
  \E^* \|\mbf - \mbZ\widehat\beta\|^2   &\leq \|\mbf-\mbZ\ti\beta\|^2 +
   \sum_{j=1}^{\ti K -K}\frac{\|\ti\beta\|_1^2}{(\ti K-j)^{2}}
   \\
   &\leq \|\mbf-\mbZ\ti\beta\|^2 + \|\ti\beta\|_1^2
   \left(\frac{1}{K-1} - \frac{1}{\ti K-1}\right) \, .
 }
This proves part (i) of the theorem. Part (ii) follows easily from
the definition of our procedure, since $p_j=1$ for all the indices
$j$ corresponding to $\ti\beta_{(1)},\dots,\ti\beta_{(k)}$ and the
$\ell_1$ norm of the estimator is preserved on every step of the
iterations. The same argument holds for part (iii) of the theorem.
\end{proof}

\medskip

Consider now the linear model of Section \ref{sec:exam}. Let $\ti
\beta$ be an estimator of parameter $\beta_0$. Using Theorem
\ref{th1} with $\mbf= \mbZ \beta_0$ we get the following bound on
the prediction loss of the randomized estimator $\widehat\beta$:
\eqsplit[pred]{\E^*\|\mbZ(\widehat\beta-\beta_0)\|^2 \le
\|\mbZ(\ti\beta-\beta_0)\|^2+\|\ti\beta\|_1^2
   \left(\frac{1}{K-1} - \frac{1}{\ti K-1}\right)\, .}
We see that if $K$ is large enough and the norm $\|\ti\beta\|_1^2$
is bounded, the difference between the losses of $\ti\beta$ and
$\widehat\beta$ is on the average not too large. For
$\ti\beta=\hat\beta_{L}$ we can replace $\|\ti\beta\|_1^2$ by $T^2$
in \eqref{pred}.

As $\ti\beta$ we may also consider another LASSO type estimator
which is somewhat different from $\hat\beta_{L}$ described in
Section \ref{sec:exam}:
\eqsplit[las2]{\ti\beta = \mathop{\arg
\min}_{\beta\in\scr^p} \left\{ \|\mby-\mbZ\beta\|^2 + r \|\beta\|_1
\right\},
}
 where $r=A \sqrt{(\log p)/n}$ with some constant $A>0$
large enough. As shown in \cite{brt07}, for this estimator, as well
as for the associated Dantzig selector, under general conditions on
the design matrix $\mbZ$ the $\ell_1$ norm satisfies
$\|\ti\beta\|_1^2= \|\beta_0\|_1^2 + o_p(s\sqrt{(\log p)/n})$ where
$s$ is the number of non-zero components of $\beta_0$. Thus, if
$\beta_0$ is sparse and has a moderate $\ell_1$ norm, the bound
\eqref{pred} can be rather accurate.

Furthermore, Theorem \ref{th1} can be readily applied to
nonparametric regression model
$$\mby=\mbf+\varepsilon$$
where $\mbf= (f(\mbX_1),\dots,f(\mbX_n))'$ and $f$ is an unknown
regression function. In this case $\mbZ\beta =f_\beta(\mbX)$ is an
approximation of $f(\mbX)$, for example as the one discussed in the
Introduction. Then, taking as $\ti\beta$ either the LASSO estimator
\eqref{las2} or the associated Dantzig selector we get immediately
sparsity oracle inequalities for prediction loss of the
corresponding randomized estimator $\widehat\beta$ that mimic (to
within the residual term $O(\|\ti\beta\|_1^2/K)$) those obtained for
the LASSO in \cite{btw07,brt07} and for the Dantzig selector in
\cite{brt07}.

It is interesting to compare our procedure with the randomization
device usually referred to as the ``Maurey argument". It is
implemented as a tool to prove approximation results over convex
classes of functions \cite{b93}. Maurey's randomization has been
used in statistics in connection to convex aggregation \cite{n00},
pages 192--193 ($K$-concentrated aggregation), and \cite{btw07},
Lemma B.1.

The Maurey randomization can be also applied to our setting. Define
the estimator $\widehat\beta_M$ as follows:
\begin{itemize}
\item[(i)] choose $K<\ti K$; draw independently at
random $K$ coordinates from $\sci$ with the probability distribution
$\{|\ti\beta_i|/\|\ti\beta\|_1\}_{i\in\sci}$,
\item[(ii)] set the $j$th coordinate of $\widehat\beta_M$ equal to
\begin{eqnarray*} \widehat\beta_{Mj} = \begin{cases}
\|\ti\beta\|_1 k_j/K & \text{ if }
\ti\beta_j> 0,\\
-\|\ti\beta\|_1 k_j/K & \text{ if }
\ti\beta_j< 0,\\
0  & \text{ if } j\not\in \sci
\end{cases}\end{eqnarray*}
where $k_j\le K$ is the number of times the $j$th coordinate is
selected at step (i).
\end{itemize}
Note that, in general, none of the non-zero coordinates of
$\widehat\beta_M$ is equal to the corresponding coordinate of the
initial estimator $\ti\beta$. The prediction risk of
$\widehat\beta_M$ is on the average not too far from that of
$\ti\beta$ as the next theorem states.
\begin{theorem}\label{th2} Under the assumptions of Theorem
\ref{th1} the randomized estimator $\widehat\beta_M$ with at most
$K<\ti K$ non-zero coordinates satisfies
\eqsplit[maurey]{\E^*\|\mbf -\mbZ\widehat\beta_M\|^2 \le \|\mbf
-\mbZ \ti\beta\|^2 +
   \frac{\|\ti\beta\|_1^2}{K}\, .}
\end{theorem}
\begin{proof} Let $\eta_1,\dots,\eta_K$ be i.i.d. random variables
taking values in $\sci$ with the probability distribution
$\{|\ti\beta_i|/\|\ti\beta\|_1\}_{i\in\sci}$. We have
$k_j=\sum_{s=1}^K I(\eta_s=j)$ where $I(\cdot)$ is the indicator
function. It is easy to see that $\E^*(\widehat\beta_{Mj}) =
\beta_j$ and the randomization covariance matrix
$\Sig^*=\E^*[(\widehat\beta_{M}-\ti\beta)(\widehat\beta_{M}-
\ti\beta)']$ has the form
\eqsplit[ranvar1]{\Sig^*=\frac{\|\ti\beta\|_1}{K}{\rm diag}
|\ti\beta_i|-\frac{1}{K}|\ti\beta| |\ti\beta|' }
where $|\ti\beta|$ is the vector of absolute values $|\ti\beta_i|$.
Acting as in \eqref{ranvar} and using \eqref{ranvar1} we get
\eqsplit{
    \E^* \|{\bf f} - \mbZ\beta_M\|^2
    &= \|{\bf f}-\mbZ\ti\beta\|^2 + {1\over
    n}\sum_{i=1}^n\mbz_i'\Sig^*\mbz_i\\
    &\leq \|{\bf f}-\mbZ\ti\beta\|^2 +
    \frac{\|\ti\beta\|_1}{K}\sum_{j\in \sci}|\ti\beta_j|
  }
which yields the result.
\end{proof}

\medskip

The residual term in \eqref{maurey} is of the same order of
magnitude $O(\|\ti\beta\|_1^2/K)$ as the one that we obtained in
Theorem \ref{th1}. In summary, $\widehat\beta_M$ does achieve the
properties (A) and (B) mentioned at the beginning of this section,
but not the property (C): it does not preserve the largest
coefficients of $\ti\beta$.

Finally, note that applying \eqref{changedLS} with $\mbf=\mby$ we
get an inequality that links the residual sums of squares (RSS) of
$\beta^*$ and $\ti\beta$:
\eqsplit[rss]{\E^*\|\mby - \mbZ\beta^*\|^2 \le \|\mby -
\mbZ\ti\beta\|^2+
   \frac{\|\ti\beta\|_1^2}{(\ti K-1)^2}\, .
   }
The left hand side of \eqref{rss} is bounded from below by the
minimum of the RSS over all the vectors $\beta$ with exactly $\ti
K-1$ non-zero entries among the $\ti K$ possible positions where the
entries of the initial estimator $\ti \beta$ are non-zero. Hence,
the minimizer $\beta^{**}$ of the residual sums of squares $\|\mby -
\mbZ\beta\|^2$ over all such $\beta$ is an estimator whose RSS does
not exceed the right hand side of \eqref{rss}. Note that
$\beta^{**}$ is obtained from $\ti\beta$ by dropping the coordinate
which has the smallest contribution to $R^2$. Iterating such a
procedure $\ti K-K$ times we get nothing but a standard backward
selection. This is exactly what we apply in Section
\ref{sec:examples}. However, the estimator obtained by this
non-randomized procedure has neither of the properties stated in
Theorem \ref{th1}since we have only a control of the RSS but not
necessarily of the prediction loss, and the $\ell_1$ norm of the
estimators is not preserved from step to step, on the difference
from our randomized procedure.


\section{Examples} \label{sec:examples}

We consider here two examples of application of our method. The
first one deals with simulated data.

\begin{example}\label{ex:simulation}
We considered a sample of size 250 from $(Y,X_1,\dots,X_{10})$,
where $X_1,\dots,X_{10}$ are \iid standard uniform,
$Y=\beta_1\ind(\frac18<X_1\leq \frac14
)+\beta_2\ind(\frac18<X_2\leq\frac12)\ind(\frac18<X_3\leq\frac38)\ind(\frac18\leq
X_4\leq \frac58)+\eps$, where $\ind(\cdot)$ denotes the indicator
function and $\eps$ is normal with mean 0 and variance such that the
population $R^2$ is 0.9. The coefficients $\beta_1$ and $\beta_2$
were selected so that the standard deviation of the second term was
three times that of the first.

We followed the hierarchical method (i)--(iii) of the Introduction.
Our initial set $\scf_0$ was a collection of $L=32$ step functions
for each of the ten variables ($d=10$). The jump points of the step
functions were equally spaced on the unit interval. The cardinality
of $\scf_0$ was 279 (after taking care of multicolinearity). At each
step we run the LASSO path until $\ti K =40$ variables were
selected, from which we selected $K=20$ variables by the standard
backward procedure. Then the model was enlarged by including
interaction terms, and the iterations were continued until there was
no increase in $R^2$.

The first step (with single effects only) ended with $R^2=0.4678$,
and the correlation of the predicted value of $Y$ with the true one
was 0.4885. The second iteration (two way interactions) ended with
$R^2=0.6303$ and correlation with the truth of 0.6115. The third
(three and four ways interactions were added) ended with
$R^2=0.7166$ and correlation of 0.5234 with the truth.  The process
stopped after the fifth step. The final predictor had correlation of
0.5300 with the true predictor.

The LASSO regularization path for the final (fifth) iteration is
presented in Figure \ref{Fig:Simu1}. The list of 20 terms included
in the model is given in the legend where $i_k$ denotes the the
$k$th step function of variable $i$. The operator $\times$ denotes
interaction of variables. We can observe that the first 12 selected
terms are functions of variables 1 to 4 that are in the true model.
Some of the 20 terms depend also on two other variables (8 and 10)
that do not belong to the true model.



    \begin{figure}[t]%
        \begin{center}%
\mbox{\makebox[\fpColWidth]{\includegraphics[width=1.6\sglWidth]
{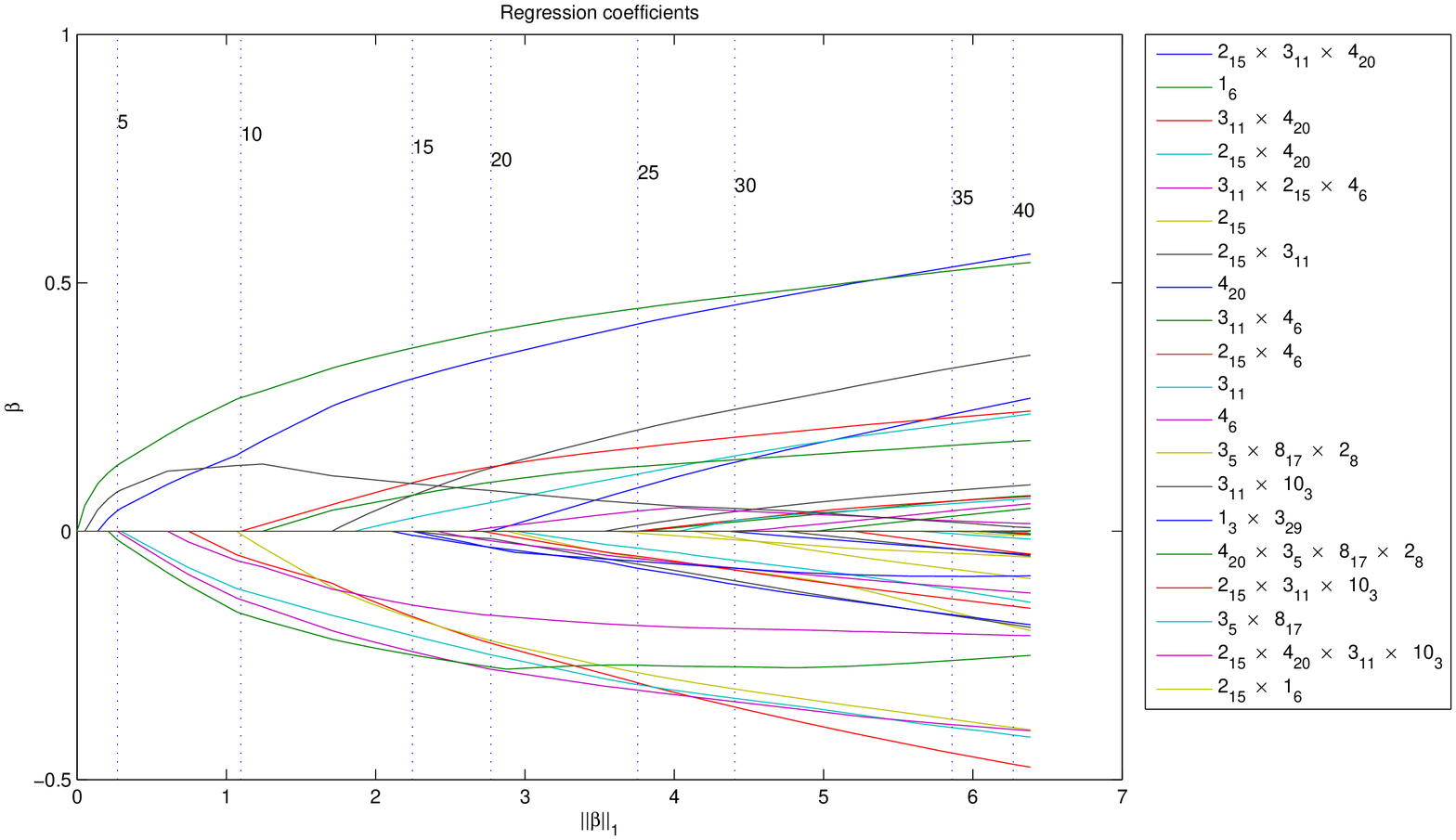} }}%
            \end{center}%
\caption{The final path of the LASSO algorithm for the simulation of Example \ref{ex:simulation}.}%
        \label{Fig:Simu1}%
    \end{figure}%

%

\end{example}

    \begin{figure}[t]%
        \begin{center}%
        \mbox{\makebox[\fpColWidth]{\includegraphics[width=1.6\sglWidth]{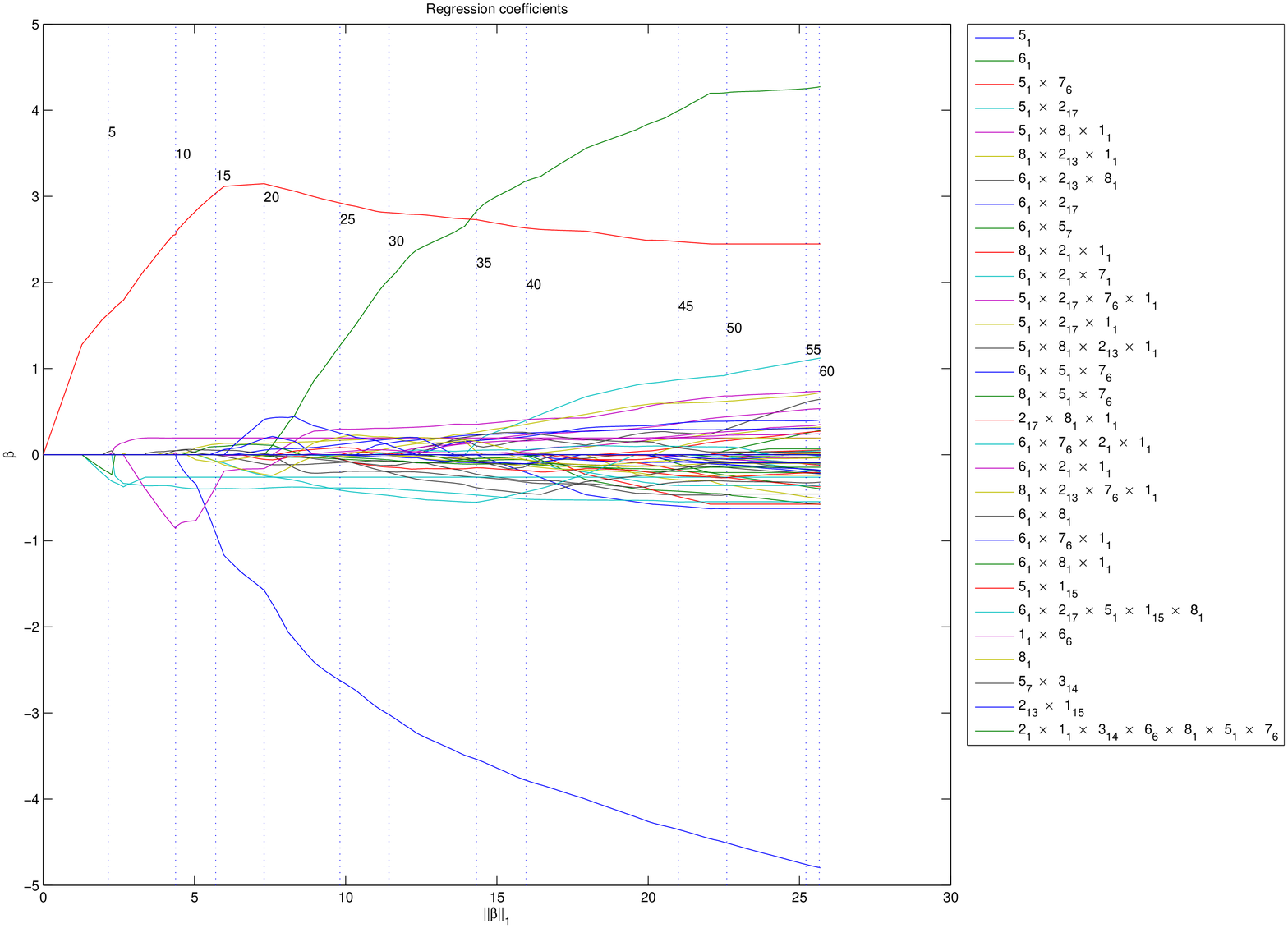} }}%
            \end{center}%
        \caption{The final path of the LASSO algorithm for the abalone data set.}%
        \label{Fig:Abalone1}%
    \end{figure}%

\begin{example}[The Abalone Data]
The abalone data set, taken from
{\tt ftp://ftp.ics.uci.edu/pub/machine-learning-databases/abalone/},
\newline gives the age of abalone (as determined by cutting
the shell and counting the number of rings) and some physical
measurements (sex, length, diameter, height, whole weight, weight of
meat, gut weight, and shell weight after being dried). The data was
described initially by Nash, et al in 1994. We selected at random
3500 data points as a training set. The 677 remaining points were
left as a test bed for cross-validation.

We used as a basic function of the univariate variable the ramp
function $(x-a)\ind(x>a)$. The range of the variables was initially
normalized to the unit interval, and we considered all break points
$a$ on the grid with spacing $1/32$. However, after dropping all
transformed variables which are in the linear span of those already
found, we were left with only 17 variables. We applied the procedure
with LASSO which ends with at most $\ti K =60$ variables, from which
at most $K=30$ were selected by backward regression.

The first stage of the algorithm ends with $R^2= 0.5586$ (since we
started with 17 terms and we were ready to leave up to 30 terms,
nothing was gained in this stage). The second stage, with all
possible main effects and  two-way interactions, dealt already with
70 variables and finished with only slightly higher $R^2$ (0.5968).
The algorithm stopped after the fifth iteration. This iteration
started with 2670 terms, and ended with $R^2=0.5779$. The
correlation of the prediction with the observed age of the test
sample was 0.5051. The result of the last stage is given in Figure
\ref{Fig:Abalone1}. It can be seen that the term with the largest
coefficient is that of the whole weight. Then come 3 terms involving
the meat weight, and its interaction with the length. The shell
weight which was most important when no interaction terms were
allowed, became not important when the interactions were added.

\end{example}

\bibliographystyle{newapajr}

\bibliography{BRT_SMR6}

\end{document}